\documentclass{article}

\usepackage[T1]{fontenc}
\usepackage{amsmath}
\usepackage{amssymb}
\usepackage{fancyhdr}

\begin{document}

\begin{center}
\textbf{ALGEBRAIC SYSTEMS OF MATRICES AND GRÖBNER BASIS}\\
\vspace{0cm}\hspace{0cm}
\end{center}

\begin{center}
\footnotesize{GERALD BOURGEOIS August 21,2007}
\end{center}

\lfoot{\vspace{-.7cm}\footnotesize\textit{Key words}: systems of
matrices, Gr\"obner basis, quaternions, nilpotent semi-group.}
\thispagestyle{fancy}

\begin{center}
\makebox[10.35cm][s]{\small\noindent\textsc{Abstract}. One studies
a particular algebraic system where the unknowns}\\
\makebox[10.35cm][s]{\small are matrices. We solve this system
according to the parameters values}\\
\makebox[5.35cm][s]{\small thanks to the theory of
Gröbner basis.                                                                             }\\
\end{center}
\vspace{0cm}\hspace{0cm}

\begin{center}
{\large 1. Introduction}
\end{center}

\noindent\textbf{1.1 The problem in question.}  $\mathbb{K}$ is a
commutative field of characteristic $0$ and algebraically closed; $n\ge2$; let $I_n$ and $0_n$ be the identity and null $(n,n)$ matrices.\\
\indent Some $(n,n)$ matrices are said $ST$ iff they are simultaneously similar to upper-triangular matrices.\\
\indent The parameters $\alpha,\beta,\gamma\in\mathbb{K}$ are given; one studies the symmetrical system in the $(n,n)$ matrices $a,b,c$ with coefficients in $\mathbb{K}$:\\
\indent $$\mathbb{S}:\left\lbrace\begin{array}{l}
a+b+c=\alpha{I_n}\\
a^2+b^2+c^2=\beta{I_n}\\
a^3+b^3+c^3=\gamma{I_n}
\end{array}\right.$$
\indent If $a,b,c$ constitute a solution of $\mathbb{S}$ and pairwise commutate, then it's easy to prove that $r(x)=6x^3-6\alpha{x^2}+(3\alpha^2-3\beta)x+3\alpha\beta-2\gamma-\alpha^3$ is a zero polynomial for $a,b,c$.\\
\indent The $r$'s discriminant is (up to the multiplicative factor $216$):\\
\ $dis(r)=9{\alpha^4}\beta-8{\alpha^3}\gamma-21\alpha^2\beta^2+36\alpha\beta\gamma-18\gamma^2-\alpha^6+3\beta^3$.\\
\indent $dis(r)=0\Longleftrightarrow{r}$ has a multiple root.\\
\indent One carries out in $\mathbb{S}$ the changes of functions:
$a=a_1+\dfrac{\alpha}{3}I_n, b=b_1+\dfrac{\alpha}{3}I_n,
\\c=c_1+\dfrac{\alpha}{3}I_n $. One deduces that: $
a_1+b_1+c_1=0$,
$a_1^2+b_1^2+c_1^2=(\beta-\dfrac{\alpha^2}{3}){I_n}$,
$a_1^3+b_1^3+c_1^3=\dfrac{\delta}{9}{I_n}$ where
$\delta=2\alpha^3-9\alpha\beta+9\gamma$.\\
 \textbf{Remark 1.1.1.} $\delta=0\Longleftrightarrow$ one of the roots of $r$ is the half sum of the 2 others.\\
 \indent One can thus rewrite the system as follows:\\
  $$\mathbb{S}:\left\lbrace\begin{array}{l}
a+b+c=0\\
a^2+b^2+c^2=\sigma{I_n}\\
a^3+b^3+c^3=\tau{I_n}
\end{array}\right.$$ with $r(x)=6x^3-3\sigma{x}-2\tau,
dis(r)=3(\sigma^3-6\tau^2)$ and $\delta=9\tau$.\\
 \textbf{1.2. The chosen method}. In the literature we didn't
 find any result about exact solutions of algebraic
 system where the unknowns are matrices; on the other hand there
 exist methods of numerical calculation which gives approximations
 of such solutions.\\
\indent Here first one seeks a Gröbner basis of the ideal generated by $\mathbb{S}$ in the ( $non$  $commutative$) ring of the polynomials in the unknowns $a,b,c$. For that, one uses the formal computation software Bergman ($cf$. [1]); it will be noted that this one works only on $homogeneous$ $polynomials$.\\
 \indent Initially one highlights 4 cases according to $dis(r)$ or $\delta$ is zero or not. Then one calculates a Gröbner basis in each case. In the continuation, these cases are classified from the most simple to the most complex.\\
\textbf{Remark 1.2.1.}  The non commutative ideals don't admit
necessarily finite Gröbner basis but the ideal generated by
$\mathbb{S}$ has always a finite basis as it would be also the
case, for example, for the nonhomogeneous system \\    $\{ ab=c,
bc=a,
ca=b\}$; on the other hand its associated homogeneous system $\{ ab=ct, bc=at, ca=bt, at=ta, bt=tb, ct=tc\}$ hasn't any finite Gröbner basis.\\
\textbf{1.3. Solutions of the system}. There exist 4
cases:\\
1.3.1. $The$ $generic$ $case$ $\delta{dis(r)}\neq{0}$: the
evaluations of $r(x)$ in $a,b,c$ are 0 and $a,b,c$ are
simultaneously diagonalizable. The resolution of $\mathbb{S}$ is thus brought back to the case $n=1$.\\
1.3.2. $\delta\neq{0}$ and $dis(r)=0$: the evaluations of $r(x)$
in $a,b,c$ are 0; $a,b,c$ pairwise commutate and generally aren't
diagonalizable; one obtains the general solution $\mathbb{S}$ in
\textbf{3}.\\
1.3.3. $\delta=0$ and $dis(r)\neq{0}$: the evaluations of
$r(x)(x^2-\dfrac{\sigma}{3})$ in $a,b,c$ are 0; $a,b,c$ are
diagonalizable but generally don't commutate; one obtains the
general solution of $\mathbb{S}$ in \textbf{4}.\\
 1.3.4. $The$ $"nilpotent"$ $case$ $\delta=dis(r)=0$: the
 evaluations  of any polynomial of degree 5 in $a,b,c$ are null;
 in general $a,b,c$ don't commutate but are $ST$.\\
 \indent In view of explicit solutions one will use the J.C. Faugere's software
 FGb ($cf$. [2]) which calculates Gröbner basis of ideals in a ($commutative$)
 ring of polynomials in the unknowns in $\mathbb{K}$.\\
 \indent There exist solutions such that $x^4$ isn't a zero
 polynomial for $a,b,c$ iff $n\geq9$.\\
 \indent The algebra $<a,b,c>$, generated by $a,b,c$, is a
 nilpotent semi group; if $a^4\neq{0}$ then its nilpotency
 class is 5; one exposes some results about the associated flag: in
 particular if $n=9$, there exists, up to isomorphism, only one
 flag of length 5. The algebras $<a,b,c>$ are classified up to
 isomorphism.\\
\textbf{1.4. Special fields}. In \textbf{6} one proves that if
$\alpha,\beta,\gamma$ are real, if $n=3$ and if $r$ has only one
real root then $\mathbb{S}$ hasn't any real solution .\\
\indent Then one studies the system
$$\left\lbrace\begin{array}{l}a+b+c=0\\a^2+b^2+c^2=v\\a^3+b^3+c^3=1\end{array}\right.$$
where $v$ is a known quaternion which isn't real, and $a,b,c$ are
unknown quaternions. One obtains a necessary and sufficient
condition on $v$ so that there exist solutions which don't
commutate with $v$. To do that one uses software $SALSA$ ( $cf$.
4); this software, which is a prolongation of the FGb software,
works on the real solutions of algebraic systems and its results
are certified. \\
\textbf{1.5. Generalization}. In \textbf{7} one expounds some
ideas about the system in the matrices $(n,n)$ $a,b,c$: $\Sigma$
    $\{a^k+b^k+c^k=\alpha_k{I_n},$ $k=1,2,3,4\}$.\\
    \indent It will be shown that the solutions of the system
    $\Sigma$ have a more complex structure than those of
    $\mathbb{S}$.
\begin{center}
{\large $2$. The generic case $\delta\neq{0}$ and $dis(r)\neq{0}$}
\end{center}

\noindent\textbf{2.1 Gröbner basis}. Let $u,v$ be ($n,n$) matrices
such that
$u^2v^3=v^3u^2$.\\
Let the system
$$\mathbb{S}:\left\lbrace\begin{array}{l}
a+b+c=0\\
a^2+b^2+c^2=u^2\\
a^3+b^3+c^3=v^3
\end{array}\right.$$
$1^{st}$ $case$: $v$ is invertible and $a,b,c$ commutate with
$u^2,v^3$.\\
The Bergman software provides a Gröbner basis from which one
extracts the relations: $2.1.0$:  $2b^2+2ab+2a^2-u^2=0$, $2.1.1$:
$ab-ba=0$,\\  $2.1.2$:   $6a^3-3au^2-2v^3=0$. \\
\textbf{Remark 2.1.1}. Let $v=I_n$; then for each $a$'s
eigenvalues $\lambda$, there exists an $u$'s eigenvalue $\mu$ such
that $6\lambda^3-3\lambda\mu^2-2=0$; it admits a multiple root if
$\mu^6=6$.\\
$2^{nd}$ $case$: $a,b,c$ commutate only with $v^3$.\\
This time we obtain relations of which those which follow:\\
$2.1.3$:   $2(a^2+b^2)+ab+ba=u^2$, $2.1.4$:   $b^3a-ab^3=ba^3-a^3b$,\\
$2.1.5$:
  $ba^3-a^3b=a^2b^2+abab-b^2a^2-baba$.\\
\noindent\textbf{2.2 Proposition 1.} Let the system in $a,b,c$
$$\mathbb{S}_1^{\ast}:\left\lbrace\begin{array}{l}
a+b+c=0\\
a^2+b^2+c^2=u^2\\
a^3+b^3+c^3=I_n
\end{array}\right.$$
where $u^2$ is a known $(n,n)$ matrix which isn't an homothety.\\
i) There always exists a solution which doesn't commutate with
$u^2$. \\
ii) If $u^2$ isn't diagonalizable, then there doesn't exist
necessarily some solution $(a,b,c)$ of $\mathbb{S}_1^{\ast}$ such
that $u$ commutates with $a,b,c$. For $n=2$ a necessary and
sufficient condition of existence of such a solution is that
$6I_2-u^6$ isn't
nilpotent of order 2.\\
$Proof$. i) It's enough to prove it for $n=2$; one is reduced to
the two following cases:\\
\[u^2=\begin{pmatrix}v&0\\0&w\end{pmatrix} \textnormal{ with } v\neq{w} \textnormal{ or } u^2=\begin{pmatrix}v&1\\0&v\end{pmatrix}\textnormal{ with }
v\neq{0}\].\\
The FGb software says that $\mathbb{S}_1^{\ast}$ always has
solutions satisfying $u^2a\neq{au^2}$. However these solutions
don't seem structured.\\
ii) $1^{st}$ case: If $u^2$ is diagonalizable then one can find a
solution such that $a,b,c$ are simultaneously diagonalizable with
$u^2$.\\
\indent The matricial polynomial $6X^3-3Xu^2-2$ has as
(matricial) resolvant polynomial $q(X)=X^2-9X+\dfrac{27}{8}u^6$.\\
In the continuation we assume $n=2$.\\
$2^{nd}$ case: Let $u$ be such that $q$ isn't a zero polynomial
for a matrix which commutates with $u$, $ie$ such that the $q$'s
matricial discriminant: $6I_2-u^6$ isn't a square. Then $6I_2-u^6$
is nilpotent of order 2 and one can be reduced to
\[u=\begin{pmatrix}e&1\\0&e\end{pmatrix} \textnormal{ where }
e^6=6; \textnormal{ let }a,b,c \textnormal{ be a solution which
commutates with }u;\]
\[a \textnormal{ is in the form: }
\begin{pmatrix}k&l\\0&k\end{pmatrix}\textnormal{ and satisfies
2.1.2; }
\begin{pmatrix}\dfrac{e^4}{3}&\dfrac{2e^3}{9}\\0&\dfrac{e^4}{3}\end{pmatrix}\text{
is the unique }\] solution for $a$. Now one seeks $b$ satisfying
2.1.1
and 2.1.0. Then necessarily\\
\[\textnormal{the discriminant
}2u^2-3a^2=\begin{pmatrix}0&\dfrac{4e}{3}\\0&0\end{pmatrix}\textnormal{must
be a square what isn't}.\]\\ $3^{rd}$ case: If $u$ is such that
$u^2$ isn't diagonalizable and $I_2-\dfrac{1}{6}{u^6}$ is the
square of a matrix $w$, then necessarily
$spectrum(u)=\{\alpha,\alpha\}$ where $\alpha\neq0$, $u^6$ isn't
diagonalizable,
$spectrum(I_2-\dfrac{1}{6}{u^6})=\{1-\dfrac{\alpha^6}{6},1-\dfrac{\alpha^6}{6}\}$
with $\alpha^6\neq6$; $w$ isn't diagonalizable and can be selected
as a polynomial in $u$,
$spectrum(I_2+w)=\{1+\epsilon\sqrt(1-\dfrac{\alpha^6}{6}),1+\epsilon\sqrt(1-\dfrac{\alpha^6}{6})\}$
where $\epsilon=\pm1$; thus $I_2+w$ isn't nilpotent and is
invertible; it's the cube of an invertible matrix $t$ what one can
choose as a polynomial in $u$. Then (Cardan's formula)
$a=6^{\frac{-1}{3}}t+6^{\frac{-2}{3}}u^2t^{-1}$ is a solution of
2.1.2 (if $w$ and $t$ are correctly selected) which commutates
with $u$ and such that $spectrum(a)=\{\lambda,\lambda\}$ where
$6\lambda^3-3\lambda\alpha^2-2=0$.\\
Existence of $b$: if $2u^2-3a^2$ is nilpotent then
$2\alpha^{2}=3\lambda^2$; the two relations between $\alpha$ and
$\lambda$ imply $\alpha^6=6$, which is contradictory. From where
we deduce the existence of a solution which commutates with
$u$.$\Box$\\
\noindent\textbf{2.3. Resolution of $\mathbb{S}$.  Theorem 1.} If
$\delta\neq0$ and $dis(r)\neq0$ then any solution $(a,b,c)$ of
$\mathbb{S}$ is such that the evaluations of $r(x)$ in $a,b,c$ are
0 and $a,b,c$ are simultaneously diagonalizable.\\ $Proof:$
according to relations 2.1.1 and 2.1.2, $a,b,c$ pairwise commutate
and $r$ is a zero polynomial for $a,b,c$; $dis(r)\neq0$ allows to
conclude.$\Box$\\

\begin{center}
{\large $3$. The case $\delta\neq{0}$ and $dis(r)={0}$}
\end{center}

\noindent\textbf{3.1 Installation}. Here $\sigma^3=6\tau^2$ and
$\tau\neq0$; using a homothety on $a,b,c$ , one is reduced to
$\tau=1$ and $\sigma^3=6$; after possible multiplication of
$a,b,c$ by $j=e^{\frac{2i\pi}{3}}$ or $j^2$ one may assume
$\sigma=\sqrt[3]{6}$ and $\tau=1$.\\
One placed oneself in $\mathbb{C}$ but the reasoning is valid in
$\mathbb{K}$.\\
Thus all the cases are brought back to a case chosen in advance.\\
Here the model is the system obtained for
$\alpha=\beta=\gamma=1$ and that we note $\mathbb{S}^{*}$.\\
\noindent\textbf{3.2. Solution of $\mathbb{S}^{*}$. Theorem 2.} If
$\delta\neq0$ and $dis(r)=0$ then the resolution of $\mathbb{S}$
is brought back to that of $\mathbb{S}^{*}$; any solution of
$\mathbb{S}^{*}$ is such that $a,b,c$ are simultaneously
similar to the block-matrices:\\
\[a'=\begin{pmatrix}I_\varphi&0&0\\0&\alpha&0\\0&0&\beta\end{pmatrix},b'=\begin{pmatrix}\gamma&0&0\\0&I_\psi&0\\0&0&-\beta\end{pmatrix},c'=\begin{pmatrix}-\gamma&0&0\\0&-\alpha&0\\0&0&I_\theta\end{pmatrix},\]
where $\alpha^2=0_\psi,\beta^2=0_\theta,\gamma^2=0_\varphi$ and
$\varphi+\psi+\theta=n$.\\
$Proof$: Relations 2.1.1, 2.1.2 remain valid: thus $a,b,c$
commutate 2 by 2 and $x^3-x^2$ is a zero polynomial for $a,b,c$.
$a^2,b^2,c^2$ are three projectors whose sum is $I_n$; by a basis
change one is reduced to:\\
\[a^2=\begin{pmatrix}I_\varphi&0&0\\0&0&0\\0&0&0\end{pmatrix},b^2=\begin{pmatrix}0&0&0\\0&I_\psi&0\\0&0&0\end{pmatrix},c^2=\begin{pmatrix}0&0&0\\0&0&0\\0&0&I_\theta\end{pmatrix},\text{ where } \varphi+\psi+\theta=n.\]
Thus
$ker(a^2)=ker(b^2-I_n)\oplus{ker(c^2-I_n)}=ker(b-I_n)\oplus{ker(c-I_n)}$
and \\
$\mathbb{K}^n=ker(a-I_n)\oplus{ker(b-I_n)}\oplus{ker(c-I_n)}=E\oplus{F}\oplus{G}$;
$a,b,c$ pairwise commutate thus $E,F,G$ are stable by $a,b,c$; in
a basis adapted to this decomposition:\\
\[a=\begin{pmatrix}I_\varphi&0&0\\0&a_1&a_2\\0&a_3&a_4\end{pmatrix},b=\begin{pmatrix}b_1&0&b_2\\0&I_\psi&0\\b_3&0&b_4\end{pmatrix},c=\begin{pmatrix}c_1&c_2&0\\c_3&c_4&0\\0&0&I_\theta\end{pmatrix}; \text{but }c=I_n-a-b\]
implies that $a,b,c$ are as announced in the statement.  $\square$

\begin{center}
{\large 4. The case $\delta=0$ and $dis(r)\neq{0}$}
\end{center}
\noindent\textbf{4.1 Gröbner basis}.

\[\text{ Let the system }\mathbb{S}_2:\left\lbrace\begin{array}{l}
a+b+c=0\\
a^2+b^2+c^2=u^2\\
a^3+b^3+c^3=0
\end{array}\right.\]
One is interested in the solutions which commutate with $u$.\\
The Bergman software provides a Gröbner basis of $\mathbb{S}_2$
from which one extracts the relations: 4.1.1 $a^2b-ba^2=0$,\\
4.1.2
 $-bab-a^2b+2a^3-u^2a=0$, 4.1.3
 $6a^5-5u^2a^3+u^4a=0.$\\
\noindent\textbf{4.2 Resolution of the system. Theorem 3.} If
$\delta=0$ and $dis(r)\neq{0}$ then every solution $a,b,c$ of
$\mathbb{S}$ is such that $r(x)(x^2-\dfrac{\sigma}{3})$ is a zero
polynomial for $a,b,c$ and $a,b,c$ are diagonalizable but
generally
they don't commutate.\\
Moreover $\mathbb{K}^n=E\oplus{F}$ where $E,F$ are $a,b,c$-stable
and such that:\\
\indent i) The $E$'s dimension $m$ is even and the restrictions of
$a,b,c$ to $E$ are \[\text{simultaneously similar to: }
A=\sqrt{\dfrac{\sigma}{3}}\begin{pmatrix}1&0\\0&-1\end{pmatrix}\otimes{I_{\frac{m}{2}}},\]
\[B=\sqrt{\dfrac{\sigma}{3}}\begin{pmatrix}-\dfrac{1}{2}&\dfrac{\sqrt{3}}{2}\\\dfrac{\sqrt{3}}{2}&\dfrac{1}{2}\end{pmatrix}\otimes{I_{\frac{m}{2}}}, C=\sqrt{\dfrac{\sigma}{3}}\begin{pmatrix}-\dfrac{1}{2}&-\dfrac{\sqrt{3}}{2}\\\dfrac{\sqrt{3}}{2}&\dfrac{1}{2}\end{pmatrix}\otimes{I_{\frac{m}{2}}}.\]
Note that these matrices aren't $ST$.\\
\indent ii) $r$ is a zero polynomial for the restrictions of
$a,b,c$  to $F$ and the later are simultaneously diagonalizable.\\
$Proof$: Here $\tau=0$, $\sigma\neq0$ and $r(x)=3x(2x^2-\sigma)$.
4.1.3 and $\sigma\neq0$ imply $\mathbb{K}^n=E\oplus{F}$ where
$E=ker(a^2-\dfrac{\sigma}{3}I_n)$ and $F=ker(r(a))$; one deduces
the first claim.\\
Proof of i): by 4.1.1 $E$ is $b,c$-stable; let $u\in{E}$ be a
$b^2$'s eigenvector: $b^2(u)=\lambda{u}$ where
$\lambda\in\{\dfrac{\sigma}{3},\dfrac{\sigma}{2},0\}$.\\
$c^2(u)=\sigma{u}-x^2(u)-y^2(u)=(\dfrac{2\sigma}{3}-\lambda)u$
where
$\dfrac{2\sigma}{3}-\lambda\in\{\dfrac{\sigma}{3},\dfrac{\sigma}{2},0\}$;
this implies that $\lambda=\dfrac{\sigma}{3}$ and if $v\in{E}$
then $b^2(v)=c^2(v)=\dfrac{\sigma}{3}v$,\\
$ker(a^2-\dfrac{\sigma}{3}I_n)\subset{ker(b^2-\dfrac{\sigma}{3}I_n)}$
hence\\
4.2.1
 $ker(a^2-\dfrac{\sigma}{3}I_n)={ker(b^2-\dfrac{\sigma}{3}I_n)}={ker(c^2-\dfrac{\sigma}{3}I_n)}$.\\
 By 4.1.2.1 $a',b',c'$, the restrictions of $a,b,c$ to $E$,
 satisfy:\\
 $-b'a'b'-a'^2b'+2a'^3-\sigma{a'}=-b'a'b'-\dfrac{\sigma}{3}b'-\dfrac{\sigma}{3}a'=0$;
 thus\\
 $0=trace(-b'^2a'+\dfrac{\sigma}{3}c')=\dfrac{\sigma}{3}trace(-a'+c')$.\\
 Then $trace(a')=trace(b')=trace(c')=0$ and $m$ is even.\\
 Thanks to a homothety one may assume that $\sigma=3$.\\
 Thanks to a change of basis we may assume that
 $a'=diag(I_{\frac{m}{2}},-I_{\frac{m}{2}})$; we seek $b'$ in the
 following form:
 \[b'=\begin{pmatrix}-\frac{1}{2}I_{\frac{m}{2}}+Y&U\\V&\frac{1}{2}I_{\frac{m}{2}}+Z\end{pmatrix};\text{ the conditions } b'^2=(a'+b')^2=I_m \text{ imply }\]
$Y=Z=0_m$ and $UV=\dfrac{3}{4}I_{\frac{m}{2}}$; thanks to a basis
change matrix in the form
\[\begin{pmatrix}P&0\\0&Q\end{pmatrix}\text{ ( leaving invariant }a'\text{) we may assume that }\]
\[b'=\begin{pmatrix}-\dfrac{1}{2}I_{\frac{m}{2}}&PUQ^{-1}\\QVP^{-1}&\dfrac{1}{2}I_{\frac{m}{2}}\end{pmatrix};\text{ there exist } P,Q\text{ such that }PUQ^{-1}=\dfrac{\sqrt{3}}{2}I_{\frac{m}{2}}\]
\[\text{therefore
}b'=\begin{pmatrix}-\dfrac{1}{2}I_{\frac{m}{2}}&\dfrac{\sqrt{3}}{2}I_{\frac{m}{2}}\\\dfrac{\sqrt{3}}{2}I_{\frac{m}{2}}&\dfrac{1}{2}I_{\frac{m}{2}}\end{pmatrix}.\]
These matrices aren't $ST$ because $(a'b'-b'a')^2=-3I_m$.\\
Proof of ii): Thanks to a homothety we may assume that $\sigma=2$.\\
$F=ker(a)\oplus{ker(a^2-I_n)}=ker(a^2)\oplus{ker(a^2-I_n)}$. By
4.1.1 $F$ is $b,c$-stable.\\
Let $a",b",c"$ be the restrictions of $a,b,c$ to $F$.\\
By 4.2.1 $r(x)=6x(x^2-1)$ is a zero polynomial for $a",b",c"$ and
the later are thus diagonalizable.\\
We can conclude thanks to the fact that $a"^2,b"^2,c"^2$ are
projectors or more quickly by the search of a Gröbner basis of the
ideal generated by
the following system in $a,b,c$ where $u$ is a known invertible matrix:\\
\[\mathbb{S}_3:\left\lbrace\begin{array}{l}
a+b+c=0\\
a^2+b^2+c^2=2u^2\\
a^3-au^2=b^3-bu^2=c^3-cu^2=0
\end{array}\right.\]
We become interested in the solutions which commutate with $u$;
the Bergman package provides a basis whose one element is
$ab-ba=0$. Thus $a",b",c"$ commutate 2 by 2 and are simultaneously
diagonalizable.$\square$\\

\begin{center}
{\large 5. The case $\delta=0$ and $dis(r)=0$}
\end{center}

\noindent\textbf{5.1 Gröbner basis}. Let the system
\[\mathbb{S}_4:\left\lbrace\begin{array}{l}
a+b+c=0\\
a^2+b^2+c^2=0\\
a^3+b^3+c^3=0
\end{array}\right.\]
The Bergman package provides a sytem equivalent to $\mathbb{S}_4$:
5.1.0  $a+b+c=0$, 5.1.1  $ab+ba=-2a^2-2b^2$, 5.1.2  $a^2b=ba^2$,
5.1.3  $2a^3=a^2b+bab$,\\ 5.1.4  $a^2ba=a^3b=-\dfrac{1}{2}a^4$,
5.1.5  $a^5=0$.\\
\noindent\textbf{5.2 Resolution of $\mathbb{S}_4$}. Let
$J_k=[(J_k)_{ij}]$ be the Jordan nilpotent matrix of dimension $k$
($(J_k)_{ij}=0$ except $(J_k)_{i,i+1}=1$).\\
\noindent\textbf{Theorem 4.} If $\delta=0$ and $dis(r)=0$ then
every solution $a,b,c$ of $\mathbb{S}$ satisfies:\\
i) $a^4=b^4; abab=baba=\dfrac{5}{2}a^4$; the value of the other
monomials of degree 4 in $(a,b)$ is $-\dfrac{1}{2}a^4$.\\
ii) $x^5$ is a zero polynomial for $a,b,c$ and more generally
every monomial of degree 5 in $(a,b,c)$ is zero.\\
iii) There exists a solution such that $a,b,c$ don't commutate 2
by 2 iff $n\geq3$.\\
iv) There exists a solution such that $a^4\neq0$ iff $n\geq9$.\\
v) $<a,b,c>$, the algebra generated by $a,b,c$, is a nilpotent
semi-group and $a,b,c$ are $ST$; if $a^4\neq0$ then its nilpotency
class is 5.\\
To this semi-group we can link a flag; for $n=9$ there exists, up
to isomorphism, only one such flag of length 5.\\
\noindent\textbf{Remark 5.2.1}. It's the single case where one
doesn't provide explicitly the general solution.\\
$Proof$: i) 5.1.3 and 5.1.4 imply
$2a^4=a^3b+abab=-\dfrac{1}{2}a^4+abab$ and $abab=\dfrac{5}{2}a^4$;
thanks to the $(a\leftrightarrow{b})$ exchange:
$baba=\dfrac{5}{2}b^4$. 5.1.3 and 5.1.4 imply\\
$2a^4=a^2ba+baba=-\dfrac{1}{2}a^4+baba$ and
$baba=\dfrac{5}{2}a^4$.\\ The remaining relations can be easily
proved.\\
ii) It's sufficient to prove the result for the monomials in
$(a,b)$ because $c=-a-b$. By i) it's sufficient to prove that
$a^5=b^5=0$; the first is 5.1.5 and the
$(a\leftrightarrow{b})$ exchange provides the second.\\
iii) If $n=2$ then the solutions are in the form:
$a=\alpha{J_2},b=\beta{J_2},c=-a-b$ where
$\alpha,\beta\in\mathbb{K}$.\\
\[\text{If } n=3 \text{ then
}a=J_3,
b=\begin{pmatrix}0&x&y\\0&0&\dfrac{-2-x}{2x+1}\\0&0&0\end{pmatrix},c=-a-b\text{
( where } x\neq\dfrac{-1}{2} \text{ and }\] $x^2+x+1\neq0)$ is a
solution whose elements don't commutate, although\\ $r(x)=x^3$ is
a zero polynomial for $a,b,c$.\\
iv) $The$ $selected$ $method$: we choose a Jordan form of $a$ and
we seek the matrix $b$; the unknowns are the $n^2$ $b$'s entries.
We rewrite the equations 5.1.1 to 5.1.5 by interchanging the roles
of $a$ and $b$; we simplify the system with the linear relations:
$a^4b=a^3ba=a^2ba^2=aba^3=ba^4=0$; we seek a Gröbner basis of the
ideal generated by this system in the ($commutative$) polynomial
ring in $n^2$ unknowns in $\mathbb{K}$. To do that we use the FGb
package of J.C. Faugère which is consistent with Maple 10.\\
$a^4\neq0$ implies $n\geq5$; we review the cases from $n=5$ to
$n=9$.\\
$n=5: a=J_5$ doesn't provide a solution for $b$.\\
$n=6: a=diag(J_5,0)$ doesn't provide a solution for $b$.\\
$n=7: a=diag(J_5,J_2)$ and $a=diag(J_5,0_2)$ don't provide a solution for $b$.\\
$n=8: a=diag(J_5,J_3),a=diag(J_5,J_2,0)$ and $a=diag(J_5,0_3)$ don't provide a solution for $b$.\\
$n=9: a=diag(J_5,J_4),a=diag(J_5,J_2,0_2),a=diag(J_5,0_4)$ and\\
$a=diag(J_5,J_2,J_2)$ don't provide a solution for $b$.\\
$a=diag(J_5,J_3,0)$ is the only choice which provides solutions;
here is a solution over $\mathbb{Q}$:
\[\begin{pmatrix}0&-\dfrac{1}{2}&0&0&0&\dfrac{3}{4}&0&0&0\\0&0&\dfrac{-1}{2}&0&0&0&\dfrac{-9}{4}&0&1\\0&0&0&\dfrac{-1}{2}&0&0&0&\dfrac{3}{4}&0\\
0&0&0&0&\dfrac{-1}{2}&0&0&0&0\\0&0&0&0&0&0&0&0&0\\0&0&-1&0&0&0&\dfrac{-1}{2}&0&0\\0&0&0&3&0&0&0&\dfrac{-1}{2}&0\\0&0&0&0&-1&0&0&0&0\\0&0&0&6&0&0&0&0&0\end{pmatrix}\]

\indent The $b$'s value relies on the choice of 9 parameters in
$\mathbb{K}$ and of 2 parameters in $\mathbb{K}$*.\\
v) The first claim is a consequence of ii),i) and of the fact that
the elements of a nilpotent semi-group are $ST$.\\
Now the definitions are those of [3].\\
\indent $Definition$: a flag $\mathcal{F}$ of length $l$ is a
filter of $\mathbb{K}^n$ subspaces: $(V_i)_{0\leq{i}\leq{l}}$ such
that $\{0\}=V_0\subsetneqq...\subsetneqq{V_l}=\mathbb{K}^n$. The
$\mathcal{F}$'s signature is the sequence:
$(dim(V_{i+1}/V_i))_i$.\\
To a nilpotent semi-group $S$ of nilpotency class $l$ we can
associate a flag of length $l$ as follows: let
$\{0\}\subset[S^{l-1}(\mathbb{K}^n)]\subset...\subset[S(\mathbb{K}^n)]\subset\mathbb{K}^n$
 where $[S^{i}(\mathbb{K}^n)]$ is the vector space generated by
 the images of the products $s_1...s_i$ where $\forall{j}\leq{i}$
 $s_j\in{S}$.\\
 \indent If $a^4\neq0$ then the flag associated to $<a,b,c>$ is
 of length 5 and by i) $V_1=a^4(\mathbb{K}^n)$; moreover
 $\forall{i}\leq4$      $V_i=[a(V_{i+1},b(V_{i+1})]$.\\
 Now we look over the particular case where $n=9$ and the flag is
 of length 5:\\
 \indent We have seen in iv) that necessarily $a$ ( therefore also
 $b$ and $c$ because $a^4=b^4=c^4$) is similar to
 $diag(J_5,J_3,0)$; with Maple 10 we prove that if
 $a=diag(J_5,J_3,0)$ then the flag associated to $<a,b,c>$ doesn't
 depend on the $b$'s choice and is:
 $V_1=[e_1],V_2=[e_1,e_2,e_6],V_3=[e_1,e_2,e_3,e_6,e_7,e_9],\\V_4=[e_1,e_2,e_3,e_4,e_6,e_7,e_8,e_9]$
 where $\{e_1,...,e_9\}$ is the $\mathbb{K}^n$ canonical basis.\\
 Thus (1,2,3,2,1) is the common signature of the flags of length 5
 if $n=9$.\\
 \indent From the precedent calculation we deduce that $a$ and $b$
 become upper triangular in the basis
 $\{e_1,e_2,e_6,e_3,e_7,e_9,e_4,e_8,e_5\}$.  $  \square$\\
\noindent\textbf{Remark 5.2.2}. If $a$ and $b$ are similar, it
doesn't imply that $c$ is similar to $a$ as this counterexample
proves it: let $n=4,a=J_4$; then (up to order) $b$ is similar to
$a$, and $c$ is similar to $diag(J_2,J_2)$.
\\
\noindent\textbf{Remark 5.2.3}. The fact that some solutions
$a,b,c$ are found in the same fixed similarity class doesn't imply
that linked flags have the same signature as we can see if we
choose $n=4$: to the $diag(J_2,J_2)$ similarity class we can link
flags of signature $(2,2)$ or $(1,2,1)$.\\
\noindent\textbf{5.3. The algebra $<a,b,c>$}. It's a nilpotent
algebra over $\mathbb{K}$ whose\\
$\{a,b,ab,ba,a^2,aba,ab^2,bab,a^2b,a^4\}$ is a vectorial
generator; therefore its dimension is at most 10.\\
\indent Again $n=9$ and $a=diag(J_5,J_3,0)$: then
$\{a,b,ab,ba,a^2,aba,ab^2,a^4\}$ is a $<a,b,c>$'s basis and this
algebra is of dimension 8; moreover its center has
$\{a^2,b^2,aba,ab^2,a^4\}$ as a basis and is of dimension 5. Yet
the algebras obtained according to the choice of $b=[b_{i,j}]$
aren't isomorphic because $bab=-aba-4ab^2+\varpi{a^4}$ and
$a^2b=ab^2-\dfrac{1}{2}\varpi{a^4}$ where
$\varpi=3b_{6,8}-b_{2,4}$; these algebras are characterized, up to
isomorphism, by the $\varpi$'s values. The $b$'s value chosen in
the theorem 4 iv)'s proof
was obtained for $\varpi=0$.\\

\begin{center}
{\large 6. Special fields}
\end{center}

\noindent\textbf{6.1 The real matrices case}. Here we prove that
the
real variant of $\mathbb{S}$ may have no solution.\\
\indent\textbf{Proposition 2}. Let $\alpha,\beta,\gamma$ be known
reals; we consider the system:
\[\mathbb{S}_{\mathbb{R}}:\left\lbrace\begin{array}{l}
a+b+c=\alpha{I_n}\\
a^2+b^2+c^2=\beta{I_n}\\
a^3+b^3+c^3=\gamma{I_n}
\end{array}\right.\] where $a,b,c$ are unknown $(n,n)$ $real$
matrices.\\
Assume that the polynomial $r$ has only one real root:\\ then
$\mathbb{S}_{\mathbb{R}}$ has at least one solution iff $n$ is even.\\
$Proof$:   $r$ has one real root $u$ and two non real roots
 $v\pm{iw}$; then $dis(r)\neq0$.\\
 a) $n$ is even.
 \[\text{If } n=2 \text{ then a solution is
 }a=uI_2,b=\begin{pmatrix}v&w\\-w&v\end{pmatrix},c=\begin{pmatrix}v&-w\\w&v\end{pmatrix}.\]
 If $n=2m$ then a solution is the tensorial product by $I_m$ of
 the former expressions.\\
\noindent\textbf{Remark 6.1.1.} More generally if $n=2$ and if the
$2^{nd}$ member of $\mathbb{S}_\mathbb{R}$ is composed of 3
matrices in $U=\{\begin{pmatrix}y&z\\-z&y\end{pmatrix}; y,z\in
\mathbb{R}\}$, a field isomorphic to $\mathbb{C}$, then there
exists a solution $(a,b,c)\in{U^3}$ which is unique, up to
order.\\
b) $n$ is odd.\\
 After a possible change of the unknown
matrices $s$ in $\dfrac{1}{|(v-u)+iw|}(s-uI_3)$ we may assume that
$u=0$ and $|v+iw|=1$.
\\
$1^{st}$ case: $\delta\neq0$ that is $v(8v^2-9)\neq0$ or ( because
$|v|<1)$ $v\neq0$: \\
According to 2.theorem 1, $a,b,c$ are simultaneously
diagonalizable over $\mathbb{C}$.\\
$2^{nd}$ case: $\delta=0$ that is $0,\pm{i}$ are the $r$'s roots.
This is the case 1.3.3 where $\sigma=-2$ and $\tau=0$:\\
According to 4.theorem 3, $F=ker(r(a))$ is a $\mathbb{C}$-space
vector of odd dimension and is $a,b,c$-stable; moreover the
restrictions of $a,b,c$ to $F$ are simultaneously diagonalizable over $\mathbb{C}$.\\
We can also see $F$ as a $\mathbb{R}$-space vector of odd
dimension which is $a,b,c$ stable.\\

\indent The following lemma proves that the both cases are
impossible:\\
\noindent\textbf{Lemma 6.1.2.} Let $u$ be a real number, $v$ be a
non real complex number and $w$ be its conjugate. Let $a,b,c$ be
real $(n,n)$ matrices and $p$ be an invertible complex matrix such
that
$a=p\;diag(a_1,\cdots,a_n)\;p^{-1},b=p\;diag(b_1,\cdots,b_n)\;p^{-1},c=p\;diag(c_1,\cdots,c_n)\;p^{-1}$
where, for all $i$, $(a_i,b_i,c_i)$ is a permutation of $(u,v,w)$.
Then necessarily $n$ is even.\\
$Proof$: (due to R. Israel) Let $U=ker(a-uI_n), U_1=\{x\in{U};
bx=vx\},U_2=\{x\in{U}; bx=wx\}$; $U=U_1\oplus{U_2}$. The complex
conjugation operator leaves $U$ invariant and interchanges $U_1$
and $U_2$. Thus $U_1$ and $U_2$ must have the same dimension over
$\mathbb{R}$, and therefore also over $\mathbb{C}$. So the
dimension of $U$ over $\mathbb{C}$ is even. Similarly
$ker(b-uI_n)$ and $ker(c-uI_n)$ have even dimensions. Since
$\mathbb{C}^n=ker(a-uI_n)\oplus{ker(b-uI_n)}\oplus{ker(c-uI_n)}$,
its dimension $n$ is also even.$\;\square$\\
\noindent\textbf{6.2 The quaternionic case }. Let $\mathcal{H}$ be
the set of quaternions; $\mathcal{H}$ is a real subspace, of basis
$\{1,i,j,k\}$, of the $(2,2)$ complex matrices. We will write a
quaternion in the form: $x=x_1+x_2i+x_3j+x_4k$.\\
\indent Let $u,v,w$ be 3 non all real quaternions which commutate;
then there exists $\rho\in{S^2}$ ($\rho^2=-1$) such that  $u,v,w$
are in the form $\lambda+\mu\rho$ where
$\lambda,\mu\in\mathbb{R}$.\\
\noindent \[\text{Now we consider the system }
\mathbb{S}_{\mathcal{H}}:  \left\lbrace\begin{array}{l}
a+b+c=u\\
a^2+b^2+c^2=v\\
a^3+b^3+c^3=w
\end{array}\right.\]
where the unknown $a,b,c$ are 3 quaternions.\\ If we are
interested only by the solutions which commutate with $\rho$ then,
thanks to the change of $\rho$ into $i\in\mathbb{C}$, we return to
the system $\mathbb{S}$ with $\mathbb{K}=\mathbb{C}$ and $n=1$: if
the $r$'s roots are $(\lambda_k+i\mu_k)_{k\le3}$ then the only
solution (up to order) in $\mathcal{H}$ which commutates with
$\rho$ is $(\lambda_k+\rho\mu_k)_{k\le3}$.\\
\indent Do there exist solutions in $\mathcal{H}$ which don't
commutate with $\rho$ ?\\
\indent We consider the following example:
$\{u=0,v\in\mathcal{H},w=1\}$.\\
$1^{st}$ case: Assume that $v\in\mathbb{R}$; as seen in 2.1. the
elements of a matricial solution commutate 2 by 2; if
$(\lambda_k+i\mu_k)_{k\le3}$ are the roots of $r(x)=6x^3-3vx-2$,
then the solutions of $\mathbb{S}_{\mathcal{H}}$ are
$(\lambda_k+\rho\mu_k)_{k\le3}$ where $\rho$ is arbitrary in
$S^2$.\\
If $v<\sqrt[3]{6}$, then $r$ has 2 non real roots and
$\mathbb{S}_{\mathcal{H}}$ has an infinity of solutions.\\
$2^{nd}$ case: We assume in the continuation that:
$\{u=0,v\in\mathcal{H}\setminus\mathbb{R},w=1\}$.\\
\indent There exists $h\in\mathcal{H}^{*}$ such that
$h^{-1}vh\in\mathbb{C}\setminus\mathbb{R}$; thus we may suppose
that $v=v_1+iv_2\in\mathbb{C}\setminus\mathbb{R}$.\\
\noindent\textbf{Remark 6.2.1.} We may assume that $v_2>0$ because
if $(a,b,c)$ is a solution associated to $v$, then
$(\bar{a},\bar{b},\bar{c})$ is a solution associated to
$\bar{v}$.\\
\indent We have seen that the solutions in $\mathcal{H}$ whose the
elements commutate with $v$ are in $\mathbb{C}$. It remains to see
if there exist some solutions in
$\mathcal{H}^3\setminus\mathbb{C}^3$. We calculate a Gröbner basis
associated to the system whose the 8 real unknowns are
$(a_i)_{i\le4},(b_i)_{i\le4}$ with $a_3\neq0$ or $a_4\neq0$, and
the real parameters are $(v_1,v_2)$.\\
\noindent\textbf{Remark 6.2.2.} If we hold a solution, the changes
of $(a_3,b_3)$ into $(-a_3,-b_3)$ or of $(a_4,b_4)$ into
$(-a_4,-b_4)$ provide other solutions, $ie$ four solutions for one
$(a_1,b_1)$'s value.\\
\indent In the following we use the Salsa package; considering an
algebraic system over $\mathbb{R}$ depending on 2 parameters, it
makes it possible to find a partition of the parameters plane such
that the number of solutions of the system in each region is a
constant; this number is then calculated by Salsa. It's because
the number of parameters  is bounded by 2 in Salsa that we don't
consider a more general second member in
$\mathbb{S}_{\mathcal{H}}$.\\
\indent Three among the basis elements are composed of polynomials
in $(a_1,b_1,v_1,v_2)$; we learn about this subsystem
$\mathcal{T}$ in the unknowns $(a_1,b_1)$ with help of Salsa. The
condition upon $v$ so that $\mathcal{T}$ has some solutions is
$\Delta=3v_1^3+4v_2^6+18v_1v_2^2-18\ge0$.\\
\indent If $\Delta>0$ then there exist six solutions in
$(a_1,b_1)$ which are associated to the permutations of a solution
in $(a_1,b_1,c_1)$.\\
\indent However the fact that $\mathcal{T}$ has some solutions
doesn't imply the existence of $\mathbb{S}_{\mathcal{H}}$'
solutions which don't commutate with $v$; there exists a second
separator curve associated to other elements of the basis; the
condition to fulfil is $v_2^2-3v_1^2\le0$.\\
\indent For example, if $v_1=-4$ then $\mathbb{S}_{\mathcal{H}}$
has, up to order, 4 solutions which don't commutate iff
$l<v_2<4\sqrt{3}$ with $l\approx2.29$.\\
\indent What is the number of solutions when $v$ is upon a
separator curve ? ( Do there exist 2 or more such solutions ?).\\

\begin{center}
{\large 7. Generalization of $\mathbb{S}$ to 4 unknowns}
\end{center}

\indent We consider the system where the unknowns are the $(n,n)$
matrices $a,b,c,d$ and where
$(\alpha_k)_{k\leq4}\subset\mathbb{K}$ is known: $\sum$
$\{a^k+b^k+c^k+d^k=\alpha_{k}I_{n}$,
$k=1,2,3,4$\}.\\
\indent Now we exhibit 2 cases in which the solutions of $\Sigma$
and $\mathbb{S}$ are very different.\\
\noindent\textbf{7.1. The generic case}. Contrary to system
$\mathbb{S}$, the calculation of a Gröbner basis of $\sum$ don't
provide any easily manageable relations, and in particular any
zero polynomial for $a$. Remember that in the case of
$\mathbb{S}$,
$a,b,c$ are simultaneously diagonalizable ( $cf$ Theorem 1).\\
\indent As in \textbf{1.1.} we can associate to the commutative
case a polynomial $r$, the form of which may be turn into
$r(x)=(x^2+ux+v)(x^2-ux+w)$ where $u,v,w\in\mathbb{K}$. If
$x^2+ux+v$ ( with roots $r_1,r_2$) is a zero polynomial for $a,c$
and $a+c=-uI_n$ and if $x^2-ux+w$ ( with roots $r_3,r_4$) is a
zero polynomial for $b,d$ and $b+d=uI_n$ then $a,b,c,d$ is a
solution of $\sum$, the elements of which aren't generally $ST$.\\
\indent Indeed, if $n=2$ then we may choose:
$a=diag(r_1,r_2),c=diag(r_2,r_1),$\\
$b=p^{-1}diag(r_3,r_4)p,d=p^{-1}diag(r_4,r_3)p$ where $p$ is
invertible. In the generic case we may choose $p$ such that
$a,b,c,d$ aren't $ST$.\\
\noindent\textbf{7.2. The "pseudonilpotent" case: $\alpha_k=0$,
$k=1,2,3,4$.} Remember that, in the case of $\mathbb{S}$, $a,b,c$
are nilpotent matrices ( $cf$ Theorem 4).\\
\noindent\textbf{Proposition 3.} If $n=2$ then the $\sum$'s
solutions satisfy one of the two following patterns ( up to
order):\\
i) $a,b,c,d$ are four nilpotent matrices such that $a+b+c+d=0$.\\
ii) $a,b,c,d$ ( defined up to a multiplicative factor) satisfy:\\
7.2.1 $\{a+b+c+d=0$, $a^2=I_2$, $b^2=jI_2$ ( where
$j=e^{\frac{2i\pi}{3}}$), $c^2=j^2I_2$, $d^2=0$,
$a+jb+j^2c=0\}$.\\
\indent If $n\geq4$ then there exist solutions which don't contain
any nilpotent matrix.\\
\noindent $Proof$: $1^{st}$ case: $a$ is neither nilpotent nor
diagonalizable; then we may assume that
$a=\begin{pmatrix}1&1\\0&1\end{pmatrix}$ and FGb proves that there
don't exist any solution.\\
$2^{nd}$ case: $a$ isn't nilpotent but is diagonalizable; then we
may assume that\\ $a=\begin{pmatrix}1&0\\0&u\end{pmatrix}$ and FGb
provides solutions; each solution satisfy 7.2.1; one of which is
this one:
$a=\begin{pmatrix}1&0\\0&-1\end{pmatrix}$,$b=\begin{pmatrix}-\frac{j^2}{2}&\frac{j-1}{2}\\\frac{1-j}{2}&\frac{j^2}{2}\end{pmatrix}$,$c=\bar{b}$,$d=\begin{pmatrix}-\frac{3}{2}&\frac{3}{2}\\-\frac{3}{2}&\frac{3}{2}\end{pmatrix}$.\\
\noindent\textbf{Remark 7.2.2.} For all $n$ if $a,b,c,d$ satisfy
7.2.1 then they constitute a $\sum$'s solution where only $d$ is
nilpotent.\\
\indent With the help of a direct sum we easily deduce the last
claim.$\square$\\
 \noindent\textbf{7.3. Commentary.} The systems
$\mathbb{S}$
and $\sum$ don't hold similar complexities.\\

\begin{center}
{\large 7 bis. A simpler example}
\end{center}

 \indent With the help of the Gröbner's basis theory, we solve the following
system where $a,b$ are $(n,n)$ unknown matrices with coefficients
in $K$:
\[\mathbb{T}:\left\lbrace\begin{array}{l}
ab+ba=I_n\\
ba^2b=0\\
\end{array}\right.\]
\noindent\textbf{Proposition 4.} If $(a,b)$ is a solution of
$\mathbb{T}$ then necessarily $n$ is even ($n=2m$) and $a,b$ are
simultaneously similar to
\[\begin{pmatrix}0&I_m\\z&0\end{pmatrix},\begin{pmatrix}0&q\\I_m&0\end{pmatrix}\]
where $z,q$ are $(m,m)$ matrices such that $zq=qz=0$\\
$Proof$: $ab$ is a projector then we may assume that
$ab=diag(I_m,0_l)$ with $m+l=n$; then $trace(ab)=m=\dfrac{n}{2}$;
$n$ is even, $l=m$ and $ba=diag(0_m,I_m)$.
\[We\; suppose\; a=\begin{pmatrix}x&y\\z&t\end{pmatrix},b=\begin{pmatrix}p&q\\r&s\end{pmatrix}\]
where $x,y,z,t,p,q,r,s$ are $(m,m)$ unknown matrices.\\
\indent With the Bergman software we study the system in the
unknowns
 $x,y,z,t,p,q,$ $r,s,u$:
 $\{xp+yr=u^2,xq+ys=0,zp+tr=0,zq+ts=0,px+qz=0,py+qt=0,rx+sz=0,ry+st=u^2,xu=ux,yu=uy,zu=uz,tu=ut,pu=up,qu=uq,ru=ur,su=us\}$.
We obtain $(xp+yr)s=0$; here $xp+yr=I_m$ and $s=0$; by the same
way $x=t=p=0$. We deduce easily that $a,b$ are simultaneously
similar to
\[\begin{pmatrix}0&y\\z&0\end{pmatrix},\begin{pmatrix}0&q\\y^{-1}&0\end{pmatrix}\;with\;zq=qz=0\;and\; y\; invertible.\]
 By a change of the $Im(ab)$-basis we may assume that $y=I_m$.               $\square$\\

\begin{center}
{\large 8. Conclusion}
\end{center}

\indent Thus we have solved the system $\mathbb{S}$ in the 3 first
cases and studied, in the $4^{th}$ case, the nilpotent algebra
generated by a solution. The use of formal calculation packages
seemed to us necessary in order to obtain such results; in
particular we mention that the use of Gröbner basis in the non
commutative case is exceptional in the literature ( the Öre
algebras case is very different because the Öre polynomials
satisfy some particular commutation's relations). Unfortunately we
befall here to the limit of computer's calculation capacity: the
$4^{th}$ case with $n=10$ imply some overflows of the FGb package;
to go further it would be necessary to work with a more effective
package ( as F5 also devised by J.C. Faugère).\\
\indent It should be interesting to study in detail the system
$\sum$ ( generalization in 4 unknowns of $\mathbb{S}$); we have
seen that the systems $\mathbb{S}$ and $\sum$ have some solutions
essentially unalike; moreover $\sum$ has Gröbner basis too
complicated for easy using; yes these basis don't triangularize
$\sum$ according to the unknowns $a,b,c$.\\
\indent On the other, the problem becomes very complicated if the
$\mathbb{S}$'s second members aren't homotheties; we have pointed
this phenomenon out during the proposition 1 ( system which has
only solutions which don't commutate) and in 6.2. where we find
requirements so that there exist in the quaternions field some non
commutating solutions; eventually it should be interesting to
study the system $\mathbb{S}_{\mathcal{H}}$ in the general case.

\begin{center}
\textsc{References}
\end{center}

\begin{enumerate}
\item \footnotesize{J. Backelin and all; the Bergman package can
be downloaded from the web site:\\
\textsf{http://servus.math.su.se/bergman/}

\item J.C. Faugere, LIP6, Paris 6; the FGb software can be downloaded from the web site:\\
\textsf{www-calfor.lip6.fr/$\thicksim$ jcf/}

\item A. Kudryavtseva and V. Mazorchuk: On the semi group of square matrices./ Arxiv preprint math.GR/0510624, 2005.
\item F. Rouillier and all; the SALSA package can be downloaded from the web site:\\
 \textsf{www-calfor.lip6.fr}  }

\end{enumerate}

\noindent\footnotesize{\textsc{G. BOURGEOIS, D\'epartement de Math\'ematiques, facult\'e de Luminy,\\
163 avenue de Luminy, case 901,\\
13288 Marseille CEDEX 09, France.}\\
\textit{E-Mail address}: \textsf{bourgeoi@lumimath.univ-mrs.fr}}

\end{document}